\documentclass[a4paper, 11pt]{article}
\usepackage{amsmath}
\usepackage{amsfonts}
\usepackage{amssymb}
\usepackage[english]{babel}
\usepackage{graphicx}
\usepackage{amsthm}
\usepackage[title]{appendix}

\usepackage{longtable}
\usepackage{tabularx}
\allowdisplaybreaks

\newtheorem{theorem}{Theorem}[section]

\newtheorem{lemma}[theorem]{Lemma}

\newtheorem{construction}[theorem]{Construction}
\def\whitebox{{\hbox{\hskip 1pt
 \vrule height 6pt depth 1.5pt
 \lower 1.5pt\vbox to 7.5pt{\hrule width
    3.2pt\vfill\hrule width 3.2pt}%
 \vrule height 6pt depth 1.5pt
 \hskip 1pt } }}
\def\qed{\ifhmode\allowbreak\else\nobreak\fi\hfill\quad\nobreak
     \whitebox\medbreak}

\newcommand{\ignore}[1]{}
\usepackage[numbers,sort&compress]{natbib}

\begin{document}
\baselineskip 16pt
\title{The Hamilton-Waterloo problem on wreath product graph $C_m \wr K_{16}$}

\author{L. Wang$^{1}$ \thanks{Corresponding author.
E-mail: {\sf hmwangli@163.com}}
and H. Cao$^{2}$  \thanks{Research
supported by the National Natural Science Foundation of China
under Grant 12071226 (H. Cao) and Suqian Science and Technology Program
under Grant Z2019110 (L. Wang).} \\
\small   1   	College of Arts and Science, Suqian University,
\small Suqian  223800,  China \\
\small   2      Institute of Mathematics, Nanjing Normal University,
\small Nanjing 210023, China
}

\date{}
\maketitle
\begin{abstract}
The Hamilton-Waterloo problem is a problem of graph factorization.
The Hamilton-Waterloo problem HWP$(H;m,n;\alpha,\beta)$
asks for a $2$-factorization of $H$ containing $\alpha$ $C_m$-factors and $\beta$ $C_n$-factors.
In this paper, we almost completely solve the Hamilton-Waterloo problem on wreath product graph $C_m \wr K_{16}$ with $C_{16}$-factors and $C_m$-factors for an odd integer $m$.


\smallskip
\end{abstract}

\section{Introduction}

In this paper, every graph will be simple. We usually denote by $V(H)$ and $E(H)$ the {\it vertex-set} and the {\it edge-set} of a graph $H$, respectively.
The complete graph on $n$ vertices is denoted by $K_n$.
We denote the cycle of length $k$ by
 $C_k$  and the complete $u$-partite graph with $u$ parts of
size $g$ by $K_u[g]$.
A {\it factor} of $H$ is a subgraph of $H$ whose vertex-set coincides with $V(H)$.
If its connected components are isomorphic to $G$, we speak of a {\it
$G$-factor}.
A {\it $G$-factorization} of $H$ is a set of edge-disjoint $G$-factors of $H$
whose edge-sets partition $E(H)$.
 A $C_k$-factorization of  $H$ is a partition of
$E(H)$ into $C_k$-factors.  The existence of a
 $C_k$-factorization of $K_u[g]$ has been completely solved
\cite{ASSW,HS,L2003,PWL,R}, which we state as the following theorem.

\begin{theorem}\label{Kv}
There exists a $C_k$-factorization of $K_u[g]$ if and only if $
g(u-1 )\equiv 0\pmod 2$, $gu\equiv 0\pmod k$, $k$ is even when
$u=2$, and $(k,u,g)\not\in\{(3,3,2),(3,6,2),(3,3,6),(6,2,6)\}$.
\end{theorem}

An $r$-regular factor is called an {\it $r$-factor}.
Actually, a 2-factor is a set of vertex-disjoint cycles. Also,
a {\it $2$-factorization} of a graph $H$ is a partition of $E(H)$
into 2-factors.
The Hamilton-Waterloo problem HWP$(H;m,n;\alpha,\beta)$
asks for a $2$-factorization of $H$ containing $\alpha$ $C_m$-factors and $\beta$
$C_n$-factors.
Let $K_v^*$ denote the complete graph $K_v$ if $v$ is odd and $K_v$ minus a
$1$-factor if $v$ is even,
we denote a solution to HWP$(K_v^*;m,n;\alpha,\beta)$ by HW$(v;m,n;\alpha,\beta)$.
Also, we denote by HWP$(H;m,n)$ the set of $(\alpha,\beta)$ for which an
HW$(v;m,n;\alpha,\beta)$ exists.
From the definition of the Hamilton-Waterloo problem,
it is not difficult to get the necessary conditions for the existence
of an HW$(v;m,n;\alpha,\beta)$ are $m | v$ when $ \alpha >0$, $n |
v$ when $ \beta >0$ and $\alpha+\beta=\lfloor \frac{v-1}{2}
\rfloor$. By Theorem~\ref{Kv}, the existence of an
HW$(v;m,n;\alpha,\beta)$ has been completely solved in the case $\alpha\beta=0$.

Most of the known results on the Hamilton-Waterloo problem for small values of $m$ and $n$.
A complete solution for the existence of an
HW$(v;3,n;\alpha,\beta)$ in the cases $n\in\{4,5,7\}$ is given in
\cite{ABBE,BB,DQS,LF,OO,WCC}.
The cases $(m,n)\in\{(3,15),(5,15),(4,6),(4,8),(4,16)$, $(8,16)\}$ are solved in \cite{ABBE}.
The existence of an HW$(v;4,n;\alpha,\beta)$ for odd
$n\geq3 $ has been solved except possibly when $v=8n$ and
$\alpha=2$, see \cite{KO,OO,WCC}.
A complete solution for $m=4$ and even $n \geq 4$ is given in \cite{FH}.
%It is shown in \cite{K} that the necessary conditions for the existence of an
%HW$(v;3,9;\alpha,\beta)$ are also sufficient except possibly when $\beta=1$.
The authors \cite{AKK} have constructed many infinite classes of HW$(v;3,3x;\alpha,\beta)$s.
In \cite{WC}, the present authors consider the Hamilton-Waterloo problem with $C_{8}$-factors and $C_m$-factors.
For more results on the Hamilton-Waterloo problem, the reader can see \cite{BDT2,BDT4,WLC}.

If $G$ and $H$ are graphs, the {\it wreath product graph} $G \wr H$ of $G$ and $H$
has a vertex-set $V(G) \times V(H)$ in which $(u_1,v_1)(u_2,v_2) \in E(G \wr H)$ whenever $u_1u_2 \in E(G)$
or, $u_1=u_2 $ and $v_1v_2 \in E(H)$.
For brevity, we denote $C_m \wr \overline{K_n}$ by $C_m[n]$, where $\overline{K_n}$ is the complement of $K_n$.

In this paper, we consider the existence of 
the Hamilton-Waterloo problem on wreath product graph $C_m \wr K_{16}$ with $C_{16}$-factors and $C_m$-factors for an odd integer $m$.

\begin{theorem}
\label{7-23}  For odd $m\geq 9$ and $7\leq r\leq 23$, the graph $C_m \wr K_{16}$
can be partitioned into $r$ $C_{16}$-factors, $23-r$ $C_m$-factors and a 1-factor.
\end{theorem}

\section{Factorizations of Cayley graphs }

We start with the definition of a Cayley graph.
Let $\Gamma$ be a finite additive group and let $S$ be a subset of
$\Gamma \backslash \{0\}$ closed under taking negatives. The {\it Cayley graph} over $\Gamma$
with connection set $S$, denoted by $\mathrm{Cay}(\Gamma, S)$, is the graph
with vertex-set $\Gamma$ and edge-set $E(\mathrm{Cay}(\Gamma, S))= \{(a,b)|
a,b \in \Gamma, a-b \in S\}$.
In this paper, we usually denote the vertex $(x,y)$ by $x_y$.

\begin{theorem}{\rm (\cite{WCC})}
\label{d=012} Let $m, n \geq 3$.  If $a\in Z_n$, $| \pm \{ 0,a,2a\} |=5$ and $(i,m)=1$, then there is a $C_m$-factorization of
$\mathrm{Cay}(Z_m \times Z_n,  \{ \pm i\}  \times  (\pm \{0, a,  2a\} ) )$.
\end{theorem}

\begin{theorem}{\rm (\cite{WC})}
\label{d}  Let $ m \geq 3$, let $n\geq4$ be even and let $0 < d < n$ be coprime to $n$. There exist
two $C_n$-factors which form a $C_n$-factorization of
$\mathrm{Cay}(Z_m \times Z_n,  \{ \pm 1\}  \times  \{ \pm d\} )$.
\end{theorem}

\begin{theorem}{\rm (\cite{WC})}
\label{d=0,4}  Let $n \geq 4$ be even and let $d=0$ ($m\geq 3$) or $d=n/2$ ($m\geq 4$ is even). There is
a $C_m$-factorization of $\mathrm{Cay}(Z_m \times Z_n,  \{ \pm 1\}  \times  \{ d\} )$.
\end{theorem}

\begin{theorem}{\rm (\cite{WC})}
\label{d=ab}  Let $ m \geq 3$ be odd, let $ n \geq 4$ be even, and let $a,b \in Z_n$ with $| \pm \{ a,b,a+b\} |=6$. There exist
six $C_m$-factors which form a $C_m$-factorization of
$\mathrm{Cay}(Z_m \times Z_n,  \{ \pm 1\}  \times (\pm \{ a,b,a+b\}) )$.
\end{theorem}

\begin{theorem}{\rm (\cite{WC})}
\label{d=0d}  Let $ m \geq 3$ be odd, let $ n \geq 4$ be even and let $1\leq d < n$. There exist three
 $C_m$-factors which form a $C_m$-factorization of
$\mathrm{Cay}(Z_m \times Z_n,  \{ \pm 1\}  \times \{ 0,\pm d\} )$.
\end{theorem}

\begin{lemma}
\label{d=n/4}  Let $ m \geq 3$ be odd and $n \equiv 0 \pmod{4}$. The graph $ \mathrm{Cay}(Z_m \times Z_{n},  \{ \pm 1\}  \times  \{ \pm n/4, n/2\} )$
can be decomposed into three $C_m$-factors.
\end{lemma}

\noindent {\it Proof:} Let $F_i=\{C_i^j \ | \ j=1,2\}$ and $C_i^j=((0,b_{i0}^j),(1,b_{i1}^j), \cdots, (m-1,b_{i,m-1}^j))$, $1\leq i \leq 3$, $j=1,2$, where

{\footnotesize \vspace{5pt}\noindent
\begin{tabular}{lllllll}
 \hspace{0.6cm} & $b_{10}^1=0$,   & $b_{11}^1=n/4$,  & $b_{12}^1=n/2;$  & $ b_{10}^2=n/4$,  & $b_{11}^2=0$,   & $b_{12}^2=-n/4; $ \\
 \hspace{0.6cm} & $b_{20}^1=0$,   & $b_{21}^1=-n/4$, & $b_{22}^1=n/4;$  & $ b_{20}^2=-n/4$, & $b_{21}^2=0$,   & $b_{22}^2=n/2; $ \\
 \hspace{0.6cm} & $b_{30}^1=0$,   & $b_{31}^1=n/2$,  & $b_{32}^1=-n/4;$ & $ b_{30}^2=-n/4$, & $b_{31}^2=n/4$, & $b_{32}^2=0. $
\end{tabular}}

\vspace{5pt}

For $m\geq 5$ and $3 \leq t \leq m-1$, let $b_{it}^j=b_{i,t-2}^j$.
Then $\{F_i +(0,s+t\cdot n/2) \ | \ 0 \leq s \leq n/4-1, 0 \leq t \leq 1\}$ is a $C_m$-factor.
\qed

\begin{lemma}
\label{d=41}  Let $l \geq 3$ and $ m \geq 2^{l-1}+1$ be odd. The graph $ \mathrm{Cay}(Z_m \times Z_{2^l},  \{ \pm 1\}  \times  \{ \pm 1, 2^{l-1}\} )$
can be decomposed into two $C_{2^l}$-factors and a $C_m$-factor.
\end{lemma}

\noindent {\it Proof:}
Let $a=1+i\cdot2^{l-2}$, $b=2^{l-1}+1+i\cdot2^{l-2}$, $c=2^{l-1}+i\cdot 2^{l-2}$, and $d=i\cdot 2^{l-2}$.

For $m=2^{l-1}+1$,
the $C_m$-factor is $\{( 0_0, 1_1, 2_2, \ldots, (2^{l-1})_{2^{l-1}})+(0, l)\ | \ l \in Z_{2^l}\}$.
Two $C_{2^l}$-factors are as below.

\vspace{5pt}

{\footnotesize
$\{( 0_0, (m-1)_1, 0_2, (m-1)_3, \ldots, 0_{2^{l}-2}, (m-1)_{2^{l}-1}),$
$( 0_a, 1_b, 2_a, 3_b, \ldots, (2^{l-1}-2)_a, (2^{l-1}-1)_b,  (2^{l-1})_c,$ \\
$ (2^{l-1}-1)_d, (2^{l-1}-2)_c, (2^{l-1}-3)_d, \ldots,2_c, 1_d)  \ | \ 0\leq i \leq 2^{l-1}-1 \}$,

$\{( (m-1)_0, 0_1, (m-1)_2, 0_3, \ldots, (m-1)_{2^{l}-2}, 0_{2^{l}-1}),  $
$( 0_{d}, 1_{2^{l-1}+d}, 2_{2^{l-1}-1+d}, 3_{2^{l-1}-2+d}, \ldots, (2^{l-1}-2)_{3+d}, $ \\ $(2^{l-1}-1)_{2+d}, (2^{l-1})_{1+d}, (2^{l-1}-1)_{2^{l-1}+1+d}, (2^{l-1}-2)_{2^{l-1}+2+d}, \ldots,2_{2^{l}-2+d}, 1_{2^{l}-1+d})  \ | \ 0\leq i \leq 2^{l-1}-1 \}$.}

\vspace{5pt}

For $m\geq 2^{l-1}+3$,
the $C_m$-factor is $\{( 0_0, 1_1, \ldots, (2^{l-1})_{2^{l-1}}, (2^{l-1}+1)_0, (2^{l-1}+2)_{2^{l-1}}, (2^{l-1}+3)_0, (2^{l-1}+4)_{2^{l-1}}, \cdots,$ $ (m-2)_0, (m-1)_{2^{l-1}}) +(0, l)\ | \ l \in Z_{2^l}\}$.
Two $C_{2^l}$-factors are given as follows.

\vspace{5pt}

{\footnotesize
$\{( 0_0, (m-1)_1, 0_2, (m-1)_3, \ldots, 0_{2^l-2}, (m-1)_{2^l-1}), $
$( 0_a, 1_b, 2_a, 3_b, \ldots, (2^{l-1}-2)_a, (2^{l-1}-1)_b, (2^{l-1})_c, (2^{l-1}-1)_d, (2^{l-1}-2)_c, (2^{l-1}-3)_d, \ldots,2_c$, $1_d)$,
$( (j+1)_0, j_1, (j+1)_2, j_3, \ldots, (j+1)_{2^l-2}, j_{2^l-1}) \ | \ 0\leq i \leq 2^{l-1}-1,  2^{l-1} \leq j \leq m-2  \}$,

\vspace{5pt}

$\{( (m-1)_0, 0_1, (m-1)_2, 0_3, \ldots, (m-1)_{2^l-2}, 0_{2^l-1}), $
$( 0_{d}, 1_{2^{l-1}+d}, 2_{2^{l-1}-1+d}, 3_{2^{l-1}-2+d}, \ldots, (2^{l-1}-2)_{3+d}, (2^{l-1}-1)_{2+d}, (2^{l-1})_{1+d}, (2^{l-1}-1)_{2^{l-1}+1+d}, (2^{l-1}-2)_{2^{l-1}+2+d}, \ldots,2_{2^{l}-2+d}$, $1_{2^{l}-1+d})$,
$( j_0, (j+1)_1, j_2, (j+1)_3, \ldots, j_{2^l-2}, (j+1)_{2^l-1}) \ | \ 0\leq i \leq 2^{l-1}-1,  2^{l-1} \leq j \leq m-2 \}$.}
\qed

\begin{lemma}
\label{d=4+K8}  Let $ m \geq 3$ and $l \geq 3$.
The graph $ \mathrm{Cay}(Z_m \times Z_{2^l},  \{ \pm 1\}  \times  \{ 2^{l-1}\} ) \cup   mK_{2^l}$
can be decomposed into $2^{l-1}$ $C_{2^l}$-factors and a 1-factor.
\end{lemma}

\noindent {\it Proof:}
By Theorem~\ref{Kv}, there exists a $C_{2^l}$-factorization of the graph $K_{2^{l-1}}[2]$.
Let $(e_1,e_2)=(0,1), (e_3,e_4)=(1+2^{l-1},2+2^{l-1}), (e_{2t+1},e_{2t+2})=(t,t+1+2^{l-1}), \ 2 \leq t \leq 2^{l-1}-2,$ $(e_{2^l-1},e_{2^l})=(2^{l-1}-1,2^{l-1})$.
Without loss of generality, let $\{ \{e_1,e_2\}, \{e_3,e_4\},\ldots, \{e_{2^l-1},e_{2^l}\} \}$ be the group set of $K_{2^{l-1}}[2]$.
There are $2^{l-1}-1$ $C_{2^l}$-factors of $K_{2^{l-1}}[2]$, denoted by $( b_{s1}, b_{s2}, \ldots, $ $ b_{s,2^l})$ for $1\leq s\leq 2^{l-1}-1$.

The required 1-factor is $\{((j+1)_{e_{4p-2}}, j_{e_{4p-1}}),(j_{e_{4p}}, (j+1)_{e_{4p+1}}) \ | \ 1 \leq p \leq 2^{l-2}, e_{2^l+1}=e_1 \}$.
Let $C_s=( 0_{b_{s1}}, 0_{b_{s2}}, \ldots, 0_{b_{s,2^l}})$, $1\leq s\leq 2^{l-1}-1$, and
$C_{2^{l-1}}=(0_{e_1}, 0_{e_2}, 1_{e_3}, 1_{e_4}, 0_{e_5}, 0_{e_6}, 1_{e_7},1_{e_8},$ $\ldots,0_{e_{2^l-3}}, 0_{e_{2^l-2}}, 1_{e_{2^l-1}},1_{e_{2^l}})$.
For $1\leq i \leq 2^{l-1}$, each $\{C_i+(l,0) \ | \ l\in Z_m\}$ is a $C_{2^l}$-factor.
Thus we get the conclusion.
\qed

For the following lemmas, we need two special 1-factorizations of $K_{16}$ with the vertex-set $Z_{16}$ whose 15 1-factors are listed as below.

 \vspace{5pt}

 {\footnotesize
\indent $I_1=\{(0,1),(3,6),(4,5),(7,10),(8,9),(11,14),(12,13),(15,2)\}$, \\
\indent $I_2=\{(2,3),(5,8),(6,7),(9,12),(10,11),(13,0),(14,15),(1,4)\}$, \\
\indent $I_3=\{(0,2),(6,1),(13,3),(7,9),(5,11),(15,12),(8,10),(14,4)\}$, \\
\indent $I_4=\{(4,6),(10,5),(1,7),(11,13),(9,15),(3,0),(12,14),(2,8)\}$, \\
\indent $I_5=\{(0,4),(10,1),(11,3),(9,2),(12,8),(14,5),(15,7),(13,6)\}$, \\
\indent $I_6=\{(6,10),(0,7),(1,9),(15,8),(2,14),(4,11),(5,13),(3,12)\}$, \\
\indent $I_7=\{(0,5),(3,7),(9,4),(2,10),(12,1),(15,11),(13,8),(6,14)\}$, \\
\indent $I_8=\{(14,3),(1,5),(7,2),(0,8),(10,15),(13,9),(11,6),(4,12)\}$, \\
\indent $I_9=\{(0,6),(2,4),(8,3),(15,5),(9,11),(7,13),(1,14),(10,12)\}$, \\
\indent $I_{10}=\{(12,2),(14,0),(4,15),(11,1),(5,7),(3,9),(13,10),(6,8)\}$, \\
\indent $I_{11}=\{(0,15),(9,5),(11,12),(6,2),(8,7),(1,13),(3,4),(14,10)\}$, \\
\indent $I_{12}=\{(10,9),(3,15),(5,6),(0,12),(2,1),(11,7),(13,14),(8,4)\}$, \\
\indent $I_{13}=\{(0,9),(6,12),(5,2),(11,8),(14,7),(4,13),(15,1),(3,10)\}$, \\
\indent $I_{14}=\{(9,6),(12,5),(2,11),(8,14),(7,4),(13,15),(1,3),(10,0)\}$, \\
\indent $I_{15}=\{(0,11),(13,2),(12,7),(14,9),(3,5),(10,4),(15,6),(1,8)\}$.
}

 \vspace{5pt}

 {\footnotesize
\indent $I_1'=\{(0,2),(1,3),(5,4),(6,8),(7,9),(11,10),(15,13),(12,14)\}$, \\
\indent $I_2'=\{(2,1),(3,5),(4,6),(8,7),(9,11),(10,15),(13,12),(14,0)\}$, \\
\indent $I_3'=\{(0,4),(2,3),(6,1),(5,7),(10,8),(14,13),(9,15),(11,12)\}$, \\
\indent $I_4'=\{(4,2),(3,6),(1,5),(7,10),(8,14),(13,9),(15,11),(12,0)\}$, \\
\indent $I_5'=\{(0,5),(8,1),(7,3),(9,2),(15,4),(14,10),(12,6),(11,13)\}$, \\
\indent $I_6'=\{(5,8),(1,7),(3,9),(2,15),(4,14),(10,12),(6,11),(13,0)\}$, \\
\indent $I_7'=\{(0,6),(2,8),(3,12),(9,14),(5,10),(1,15),(7,13),(4,11)\}$, \\
\indent $I_8'=\{(6,2),(8,3),(12,9),(14,5),(10,1),(15,7),(13,4),(11,0)\}$, \\
\indent $I_9'=\{(0,7),(12,2),(13,3),(14,6),(10,4),(9,1),(11,5),(15,8)\}$, \\
\indent $I_{10}'=\{(7,12),(2,13),(3,14),(6,10),(4,9),(1,11),(5,15),(8,0)\}$, \\
\indent $I_{11}'=\{(0,10),(2,11),(3,15),(7,14),(6,13),(5,9),(1,12),(4,8)\}$, \\
\indent $I_{12}'=\{(0,1),(7,2),(12,4),(10,3),(13,5),(11,14),(8,9),(15,6)\}$, \\
\indent $I_{13}'=\{(6,7),(13,8),(2,10),(0,9),(3,11),(1,4),(14,15),(5,12)\}$, \\
\indent $I_{14}'=\{(0,3),(13,1),(7,4),(14,2),(8,11),(5,6),(12,15),(9,10)\}$, \\
\indent $I_{15}'=\{(10,13),(7,11),(1,14),(8,12),(2,5),(15,0),(6,9),(3,4)\}$.
}

Note that both $I_{13} \cup I_{14}$ and $I_{2i-1}' \cup I_{2i}'$ $(1 \leq i \leq5)$ can form an 16-cycle.

For an integer $m \geq 2$, $mH$ denotes $m$ vertex-disjoint copies of a graph $H$.
For brevity, we use $mI_k ($or $ mI_k')$ to denote the graph with the vertex-set $Z_m \times Z_{16}$ and the edge-set $\{(j_a,j_b) \ | \ j \in Z_m, (a,b) \in I_k ($or $I_k'), a \neq b \}$ for $1\leq k \leq 15$.
Similarly, $mK_n$ denotes the graph with the vertex-set $Z_m \times Z_n$ and the edge-set $\{(j_a,j_b) \ | \ j \in Z_m, (a,b) \in E(K_n )\}$.

\begin{lemma}
\label{d=2-13}  Let $ m \geq 3$ and $i \in \{2,4,6\}$. There exist
two $C_{16}$-factors which form a $C_{16}$-factorization of
$ \mathrm{Cay}(Z_m \times Z_{16},  \{ \pm 1\}  \times  \{ i\} )  \cup  mI_{i-1} \cup mI_i$.
\end{lemma}

\noindent {\it Proof:}
For $i \in \{2,4,6\}$,
let $C_1=(0_0, 0_1, 1_3, 1_6,0_4, 0_5, 1_7, 1_{10}, 0_8, 0_9,1_{11},$ $ 1_{14}, 0_{12}, 0_{13}, 1_{15}, 1_{2})$,
$C_1=(0_0, 0_2, 1_6, 1_1, 0_{13}, 0_3, 1_7, 1_{9}, 0_5, 0_{11}, 1_{15}, 1_{12}, 0_{8}, 0_{10}, 1_{14},$ $ 1_{4})$,
 and
$C_1=(0_0, 0_4, 1_{10}, 1_1, 0_{11},$ $ 0_3, 1_9, 1_{2}, 0_{12}, 0_{8}, 1_{14}, 1_{5}, 0_{15}, 0_{7}, 1_{13}, 1_{6})$, respectively.
$C_2$ can be generated from $C_1$ by $(-, +i \pmod{16})$.
Each $\{C_t+(l, 0) \ | \  l\in Z_m\}$ is a $C_{16}$-factor for $t=1,2$.
Thus we get two $C_{16}$-factors which form a $C_{16}$-factorization of
$ \mathrm{Cay}(Z_m \times Z_{16},  \{ \pm 1\}  \times  \{ i\} ) \cup  mI_{i-1} \cup mI_i$.\qed

\begin{lemma}
\label{d=6-27}  Let $ m \geq 3$ and $i \in \{2,4,6\}$. There exist
two $C_{16}$-factors which form a $C_{16}$-factorization of
$ \mathrm{Cay}(Z_m \times Z_{16},  \{ \pm 1\}  \times  \{ -i\} )  \cup  mI_{6+i-1} \cup mI_{6+i}$.
\end{lemma}

\noindent {\it Proof:}
For $i \in \{2,4,6\}$,
let $C_1=(0_0, 0_5, 1_3, 1_7, 0_9, 0_4, 1_2, 1_{10}, 0_{12}, 0_1, 1_{15},$ $ 1_{11}, 0_{13}, 0_{8}, 1_{6}, 1_{14})$,
$C_1=(0_0, 0_6, 1_2, 1_4, 0_{8}, 0_3, 1_{15}, 1_{5}, 0_9, 0_{11}, 1_{7}, 1_{13}, 0_{1}, 0_{14}, 1_{10},$ $ 1_{12})$,
 and
$C_1=(0_0, 0_{15}, 1_{9}, 1_5, 0_{11},$ $ 0_{12}, 1_6, 1_{2}, 0_{8}, 0_{7}, 1_{1}, 1_{13}, 0_{3}, 0_{4}, 1_{14}, 1_{10})$, respectively.
$C_2$ can be generated from $C_1$ by $(-, -i \pmod{16})$.
Each $\{C_t+(l, 0) \ | \  l\in Z_m\}$ is a $C_{16}$-factor.
We get the conclusion. \qed

\begin{lemma}
\label{d=4-456}  Let $ m \geq 3$.
The graph $ \mathrm{Cay}(Z_m \times Z_{16},  \{ \pm 1\}  \times  \{ 8\} ) \cup m(\cup_{i=13}^{15} I_i)$
can be decomposed into two $C_{16}$-factors and a 1-factor.
\end{lemma}

\noindent {\it Proof:}
The 1-factor is $\{((j+1)_9, j_1), ((j+1)_{2}, j_{10}),((j+1)_{11}, j_{3}),((j+1)_{12}, j_4),((j+1)_{13}, j_{5}),$ $((j+1)_{14}, j_6),((j+1)_7, j_{15}),((j+1)_0, j_8) \ | \ j \in Z_m\}$.
Let $C_1=(0_0, 0_9, 0_6, 0_{12}, 0_5, 0_2, 0_{11}, 0_{8}, 0_{14}, $ $ 0_7, 0_{4}, 0_{13}, 0_{15}, 0_1, 0_{3}, 0_{10})$ and
$C_2=(0_0, $ $0_{11}, 1_3, 1_5, 0_{13}, 0_2, 1_{10}, 1_4, 0_{12}, 0_7, 1_{15}, 1_6, 0_{14}, 0_9, 1_1, 1_8)$.
Each $\{C_t+(l,0) \ | \ l\in Z_m\}$ is a $C_{16}$-factor.
Thus, we obtain the required design.
\qed

\begin{lemma}
\label{d=8-6}  Let $ m \geq 3$.
The graph $ \mathrm{Cay}(Z_m \times Z_{16},  \{ \pm 1\}  \times  \{ 8\} ) \cup mI_{11}'$
can be decomposed into a $C_{16}$-factor and a 1-factor.
\end{lemma}

\noindent {\it Proof:}
The 1-factor is $\{((j+1)_1, j_9), ((j+1)_{10}, j_2),((j+1)_3, j_{11}),((j+1)_{12}, j_4),((j+1)_{13}, j_{5}),$ $((j+1)_{6}, j_{14}),((j+1)_{15}, j_{7}),((j+1)_0, j_8) \ | \ j \in Z_m\}$.
Let
$C=(0_0, 0_{10}, 1_2, 1_{11}, 0_{3}, 0_{15}, 1_{7}, 1_{14}, 0_{6},$ $ 0_{13}, 1_{5}, 1_9, 0_{1}, 0_{12}, 1_4, 1_8)$,
then $\{C+(l,0) \ | \ l\in Z_m\}$ is a $C_{16}$-factor.
\qed

\begin{lemma}
\label{d=6+I}  Let $ m \geq 3$.
The graph $ \mathrm{Cay}(Z_m \times Z_{16},  \{ \pm 1\}  \times  \{ \pm 6\} ) \cup m(\cup_{i=12}^{15} I_i')$
can be decomposed into four $C_{16}$-factors.
\end{lemma}

\noindent {\it Proof:} For $1 \leq t \leq 4$, each $\{C_t+(l,0) \ | \ l\in Z_m\}$ is a $C_{16}$-factor, where

$C_1=(0_0, 0_{1}, 1_7, 1_2, 0_{12}, 0_4, 1_{10}, 1_3, 0_{13}, 0_5, 1_{11}, 1_{14}, 0_{8}, 0_9, 1_{15}, 1_6)$,

$C_2=(0_6, 0_{7}, 1_{13}, 1_8, 0_{2}, 0_{10}, 1_{0}, 1_9, 0_{3}, 0_{11}, 1_{1}, 1_4, 0_{14}, 0_{15}, 1_5, 1_{12})$,

$C_3=(0_0, 0_{3}, 1_{13}, 1_1, 0_{7}, 0_4, 1_{14}, 1_2, 0_{8}, 0_{11}, 1_{5}, 1_6, 0_{12}, 0_{15}, 1_9, 1_{10})$,

$C_4=(0_{10}, 0_{13}, 1_7, 1_{11}, 0_{1}, 0_{14}, 1_{8}, 1_{12}, 0_{2}, 0_5, 1_{15}, 1_0, 0_{6}, 0_9, 1_3, 1_4)$.\\
Then we get the required design.
\qed

 \section{Main Results}

In order to prove our main results, we give the following known results.

\begin{theorem}{\rm (\cite{CNT,PWL})}
\label{Cm}  There exists a $C_{m}$-factorization of $C_m[n]$ for $m \geq 3$ and $n \geq 1$
except for $(m,n)=(3,6)$ and $(m,n) \in \{(l,2) \ | \ l \geq 3 $ is odd $\}$.
\end{theorem}

\begin{theorem}{\rm (\cite{LD})}
\label{Cmn}  There exists a $C_{mn}$-factorization of $C_m[n]$ for $m \geq 3$ and $n \geq 1$.
\end{theorem}

\begin{theorem}{\rm (\cite{OO})}
\label{Cm4} The graph $C_m[4]$ can be decomposed into $\alpha$ $C_4$-factors and $4-\alpha$ $C_m$-factors for $m\geq 3$ and $\alpha \in \{0,2,4\}$.
\end{theorem}

\begin{construction} {\rm (\cite{WCC})}
\label{C-RGDD} If there exist an {\rm
HW}$( K_u[g] ; m ,n; \alpha,  \beta)$ and an {\rm
HW}$(g;m,n; \alpha', \beta' )$, then an {\rm
HW}$(gu;m,n;  \alpha + \alpha',  \beta +\beta' )$ exists.
\end{construction}

The proof of the following construction is similar to Construction 3.5 in \cite{WC}, so here we just state it and omit the proof.

\begin{construction}
\label{Cmnw} If $C_m[n]$ can be decomposed into $\alpha$ $C_n$-factors and $n-\alpha$ $C_m$-factors, there exist a $C_{nw}$-factorization of $C_n[w]$ and a $C_{m}$-factorization of $C_m[w]$,
then $C_m[nw]$ can be decomposed into $w\alpha$ $C_{nw}$-factors and $w(n-\alpha)$ $C_m$-factors.
\end{construction}

\begin{construction}{\rm (\cite{WC})}
\label{Cmn+mKn} If $(\alpha,\beta) \in $ {\rm HWP}$(C_m[n];m,n)$,
then $(\alpha, \beta+ \lfloor \frac{n-1}{2}\rfloor) \in $ {\rm HWP}$(C_m \wr K_n;m,n)$.
\end{construction}

\begin{lemma}
\label{0248}  For odd $m\geq 9$ and $r \in \{0,2,4,6,8,16\}$, $(r,16-r)\in$ {\rm HWP}$(C_m[16];16,m)$.
\end{lemma}

\noindent {\it Proof:} We consider the following four cases.

{\bf Case 1:} $r \in \{0,8,16\}$.

The graph $C_m[4]$ can be decomposed into $r/4$ $C_4$-factors and $4-r/4$ $C_m$-factors for $m\geq 3$ by Theorem~\ref{Cm4}.
From Theorem~\ref{Cmn}, the graph $C_4[4]$ can be partitioned into four $C_{16}$-factors.
From Theorem~\ref{Cm}, the graph $C_m[4]$ can be partitioned into four $C_{m}$-factors.
Applying Construction~\ref{Cmnw}, we can get the conclusion.

 \vspace{5pt}

{\bf Case 2:} $r=2$.

Two $C_{16}$-factors are given from a $C_{16}$-factorization of
$ \mathrm{Cay}(Z_m \times Z_{16},  \{ \pm 1\}  \times  \{ \pm 5\} )$ by Theorem~\ref{d}.
The required 14 $C_m$-factors can be obtained as follows.
The graph $ \mathrm{Cay}(Z_m \times Z_{16},  \{ \pm 1\}  \times  \{0,\pm 1,\pm 2 \} )$
can be decomposed into five $C_m$-factors by Theorem~\ref{d=012}. Similarly, $ \mathrm{Cay}(Z_m \times Z_{16},  \{ \pm 1\}  \times  (\pm \{ 3,6,7\}) )$
can be partitioned into six $C_m$-factors from Theorem~\ref{d=ab}, $ \mathrm{Cay}(Z_m \times Z_{16},  \{ \pm 1\}  \times  \{ \pm 4,8\} )$
can be partitioned into three $C_m$-factors from Lemma~\ref{d=n/4}.

 \vspace{5pt}

{\bf Case 3:} $r=4$.

Four $C_{16}$-factors are given from a $C_{16}$-factorization of
$ \mathrm{Cay}(Z_m \times Z_{16},  \{ \pm 1\}  \times  (\pm \{ 5,7\}) )$ by Theorem~\ref{d}.
The required 12 $C_m$-factors can be obtained as follows.
The graph $ \mathrm{Cay}(Z_m \times Z_{16},  \{ \pm 1\}  \times  \{0,\pm 6 \} )$
can be decomposed into three $C_m$-factors by Theorem~\ref{d=0d}.
Similarly, $ \mathrm{Cay}(Z_m \times Z_{16},  \{ \pm 1\}  \times  (\pm \{ 1,2,3\}) )$
can be partitioned into six $C_m$-factors from Theorem~\ref{d=ab},
$ \mathrm{Cay}(Z_m \times Z_{16},  \{ \pm 1\}  \times  \{\pm 4,8\})$
can be partitioned into three $C_m$-factors from Lemma~\ref{d=n/4}.

 \vspace{5pt}

{\bf Case 4:} $r=6$.

Four of the required $C_{16}$-factors are given from a $C_{16}$-factorization of
$ \mathrm{Cay}(Z_m \times Z_{16},  \{ \pm 1\}  \times  (\pm \{ 5,7\}) )$ by Theorem~\ref{d}.
From Lemma~\ref{d=41}, the graph $ \mathrm{Cay}(Z_m \times Z_{16},  \{ \pm 1\}  \times  \{\pm 1 ,8\} )$
can be decomposed into two $C_{16}$-factors and a $C_m$-factor.
The other required $C_m$-factors can be obtained as follows.
$ \mathrm{Cay}(Z_m \times Z_{16},  \{ \pm 1\}  \times  \{0,\pm 3 \} )$
can be decomposed into three $C_m$-factors by Theorem~\ref{d=0d}.
Similarly, $ \mathrm{Cay}(Z_m \times Z_{16},  \{ \pm 1\}  \times  (\pm \{ 2,4,6\}) )$
can be partitioned into six $C_m$-factors from Theorem~\ref{d=ab}.
\qed

Finally, we prove our main theorem.

\noindent \underline{\bf \it  Proof of Theorem~\ref{7-23} :}
Let the vertex-set be $Z_{m} \times Z_{16}$. We distinguish 12 cases as below.

\vspace{5pt}

{\bf Case 1}: $r \in \{7,9,11,13,15,23\}$.

For odd $m\geq 9$ and $r_1 \in \{0,2,4,6,8,16\}$, we have $(16-r_1, r_1) \in $ {\rm HWP}$(C_m[16];m,16)$ and $(0, 7) \in $ {\rm HWP}$(16;m,16)$ by Lemma~\ref{0248} and Theorem~\ref{Kv}, respectively.
Thus we get $(16-r_1, r_1+7) \in $ {\rm HWP}$(C_m \wr K_{16};m,16)$ by using Construction~\ref{Cmn+mKn}.

\vspace{5pt}

{\bf Case 2}: $r=8$.

$ \mathrm{Cay}(Z_m \times Z_{16},  \{ \pm 1\}  \times  \{ 8\} ) \cup  mK_{16}$
can be decomposed into eight $C_{16}$-factors and a 1-factor by Lemma~\ref{d=4+K8}.
Three $C_m$-factors come from a $C_m$-factorization of $ \mathrm{Cay}(Z_m \times Z_{16},  \{ \pm 1\}  \times  \{ 0, \pm 2\} )$ by Theorem~\ref{d=0d}.
$ \mathrm{Cay}(Z_m \times Z_{16},  \{ \pm 1\}  \times  ( \pm \{ 1,5,6\}) )$ and $ \mathrm{Cay}(Z_m \times Z_{16},  \{ \pm 1\}  \times  ( \pm \{ 3,4,7\}) )$
can be partitioned into 12 $C_m$-factors from Theorem~\ref{d=ab}.

\vspace{5pt}

{\bf Case 3}: $r=10$.

$ \mathrm{Cay}(Z_m \times Z_{16},  \{ \pm 1\}  \times  \{ \pm 6, 8\} ) \cup  m(\cup_{i=11}^{15} I_i')$
can be decomposed into five $C_{16}$-factors and a 1-factor by Lemmas~\ref{d=8-6} and~\ref{d=6+I}.
The graph $ m(I_{2j-1}' \cup I_{2j}')$ can be decomposed into a $C_{16}$-factor for $1 \leq j \leq 5$ since $I_{2j-1}' \cup I_{2j}'$ can form a 16-cycle.
In other words, $ \mathrm{Cay}(Z_m \times Z_{16},  \{ \pm 1\}  \times  \{ \pm 6, 8\} ) \cup  mK_{16}$
can be decomposed into ten $C_{16}$-factors and a 1-factor.
A $C_m$-factor comes from a $C_m$-factorization of $ \mathrm{Cay}(Z_m \times Z_{16},  \{ \pm 1\}  \times  \{ 0\} )$ by Theorem~\ref{d=0,4}.
$ \mathrm{Cay}(Z_m \times Z_{16},  \{ \pm 1\}  \times  ( \pm \{ 1,3,4\}) )$ and $ \mathrm{Cay}(Z_m \times Z_{16},  \{ \pm 1\}  \times  ( \pm \{ 2,5,7\}) )$
can be partitioned into 12 $C_m$-factors from Theorem~\ref{d=ab}.

\vspace{5pt}

{\bf Case 4}: $r=12$.

$ \mathrm{Cay}(Z_m \times Z_{16},  \{ \pm 1\}  \times  \{ \pm 6, 8\} ) \cup  mK_{16}$
can be decomposed into ten $C_{16}$-factors and a 1-factor.
 $ \mathrm{Cay}(Z_m \times Z_{16},  \{ \pm 1\} \times \{ \pm 5\} )$
can be decomposed into two $C_{16}$-factors by Theorem~\ref{d}.
Five $C_m$-factors come from a $C_m$-factorization of $ \mathrm{Cay}(Z_m \times Z_{16},  \{ \pm 1\}  \times  \{ 0, \pm 1, \pm 2 \} )$ by Theorem~\ref{d=012}.
$ \mathrm{Cay}(Z_m \times Z_{16},  \{ \pm 1\}  \times  ( \pm \{ 3,4,7 \}) )$ can be partitioned into six $C_m$-factors from Theorem~\ref{d=ab}.

\vspace{5pt}

{\bf Case 5}: $r=14$.

$ \mathrm{Cay}(Z_m \times Z_{16},  \{ \pm 1\}  \times  \{ 8\} ) \cup  mK_{16}$
can be decomposed into eight $C_{16}$-factors and a 1-factor by Lemma~\ref{d=4+K8}.
The last six $C_{16}$-factors come from a $C_{16}$-factorization of $ \mathrm{Cay}(Z_m \times Z_{16},  \{ \pm 1\}  \times  (\pm \{ 3,5,7\} ) )$ by Theorem~\ref{d}.
Three $C_m$-factors come from a $C_m$-factorization of $ \mathrm{Cay}(Z_m \times Z_{16},  \{ \pm 1\}  \times  \{ 0, \pm 1\} )$ by Theorem~\ref{d=0d}.
$ \mathrm{Cay}(Z_m \times Z_{16},  \{ \pm 1\}  \times  ( \pm \{ 2,4,6\}) )$ can be partitioned into six $C_m$-factors from Theorem~\ref{d=ab}.

\vspace{5pt}

{\bf Case 6}: $r=16$.

$ \mathrm{Cay}(Z_m \times Z_{16},  \{ \pm 1\}  \times  \{ 8\} ) \cup  mK_{16}$
can be decomposed into eight $C_{16}$-factors and a 1-factor by Lemma~\ref{d=4+K8}.
The last eight $C_{16}$-factors come from a $C_{16}$-factorization of $ \mathrm{Cay}(Z_m \times Z_{16},  \{ \pm 1\}  \times  (\pm \{ 1,3,5,7\} ) )$ by Theorem~\ref{d}.
A $C_m$-factor comes from a $C_m$-factorization of $ \mathrm{Cay}(Z_m \times Z_{16},  \{ \pm 1\}  \times  \{ 0\} )$ by Theorem~\ref{d=0,4}.
$ \mathrm{Cay}(Z_m \times Z_{16},  \{ \pm 1\}  \times  ( \pm \{ 2,4,6\}) )$ can be partitioned into six $C_m$-factors from Theorem~\ref{d=ab}.

\vspace{5pt}

{\bf Case 7}: $r=18$.

$ \mathrm{Cay}(Z_m \times Z_{16},  \{ \pm 1\}  \times  \{ \pm 6, 8\} ) \cup  mK_{16}$
can be decomposed into ten $C_{16}$-factors and a 1-factor.
$ \mathrm{Cay}(Z_m \times Z_{16},  \{ \pm 1\} \times (\pm \{ 1,3,5,7\}) )$
can be decomposed into eight $C_{16}$-factors by Theorem~\ref{d}.
The last five $C_m$-factors come from a $C_m$-factorization of $ \mathrm{Cay}(Z_m \times Z_{16},  \{ \pm 1\}  \times  \{ 0, \pm 2, \pm 4 \} )$ by Theorem~\ref{d=012}.

\vspace{5pt}

{\bf Case 8}: $r=20$.

$ \mathrm{Cay}(Z_m \times Z_{16},  \{ \pm 1\} \times (\pm \{ 3,5,7\}) )$
can be decomposed into six $C_{16}$-factors by Theorem~\ref{d}.
$ \mathrm{Cay}(Z_m \times Z_{16},  \{ \pm 1\}  \times (\pm \{ 2,4,6\} ) ) \cup  m(\cup_{i=1}^{12} I_i)$
can be decomposed into 12 $C_{16}$-factors by Lemmas~\ref{d=2-13} and~\ref{d=6-27}.
$ \mathrm{Cay}(Z_m \times Z_{16},  \{ \pm 1\}  \times  \{ 8\} ) \cup  m(\cup_{i=13}^{15} I_i)$
can be decomposed into two $C_{16}$-factors and a 1-factor by Lemma~\ref{d=4-456}.
Three $C_m$-factors come from a $C_m$-factorization of $ \mathrm{Cay}(Z_m \times Z_{16},  \{ \pm 1\}  \times  \{ 0, \pm 1\} )$ by Theorem~\ref{d=0d}.

\vspace{5pt}

{\bf Case 9}: $r=22$.

$ \mathrm{Cay}(Z_m \times Z_{16},  \{ \pm 1\} \times (\pm \{ 1,3,5,7\}) )$
can be decomposed into eight $C_{16}$-factors by Theorem~\ref{d}.
$ \mathrm{Cay}(Z_m \times Z_{16},  \{ \pm 1\}  \times (\pm \{ 2,4,6\} ) ) \cup  m(\cup_{i=1}^{12} I_i)$
can be decomposed into 12 $C_{16}$-factors by Lemmas~\ref{d=2-13} and~\ref{d=6-27}.
$ \mathrm{Cay}(Z_m \times Z_{16},  \{ \pm 1\}  \times  \{ 8\} ) \cup  m(\cup_{i=13}^{15} I_i)$
can be decomposed into two $C_{16}$-factors and a 1-factor by Lemma~\ref{d=4-456}.
A $C_m$-factor comes from a $C_m$-factorization of $ \mathrm{Cay}(Z_m \times Z_{16},  \{ \pm 1\}  \times  \{ 0\} )$ by Theorem~\ref{d=0,4}.

\vspace{5pt}

{\bf Case 10}: $r=17$.

The 1-factor is $mI_{11}'$.
$ \mathrm{Cay}(Z_m \times Z_{16},  \{ \pm 1\}  \times  \{ \pm 6\} ) \cup  m(\cup_{i=12}^{15} I_i')$
can be decomposed into four $C_{16}$-factors by Lemma~\ref{d=6+I}.
The graph $ m(I_{2j-1}' \cup I_{2j}')$ can be decomposed into a $C_{16}$-factor for $1 \leq j \leq 5$, since $I_{2j-1}' \cup I_{2j}'$ can form a 16-cycle.
$ \mathrm{Cay}(Z_m \times Z_{16},  \{ \pm 1\}  \times  ( \pm \{ 3,5,7\}) )$
can be partitioned into six $C_{16}$-factors from Theorem~\ref{d}.
$ \mathrm{Cay}(Z_m \times Z_{16},  \{ \pm 1\}  \times  \{ \pm 1, 8\} ) $
can be decomposed into two $C_{16}$-factors and a $C_{m}$-factor by Lemma~\ref{d=41}.
The last five $C_m$-factors come from a $C_m$-factorization of $ \mathrm{Cay}(Z_m \times Z_{16},  \{ \pm 1\}  \times  \{ 0, \pm 2, \pm 4\} )$ by Theorem~\ref{d=012}.

\vspace{5pt}

{\bf Case 11}: $r=19$.

The 1-factor is $mI_{15}$.
$ \mathrm{Cay}(Z_m \times Z_{16},  \{ \pm 1\}  \times (\pm \{ 2,4,6\} ) ) \cup  m(\cup_{i=1}^{12} I_i)$
can be decomposed into 12 $C_{16}$-factors by Lemmas~\ref{d=2-13} and~\ref{d=6-27}.
The graph $ m(I_{13} \cup I_{14})$ can be decomposed into a $C_{16}$-factor since $I_{13} \cup I_{14}$ can form a 16-cycle.
$ \mathrm{Cay}(Z_m \times Z_{16},  \{ \pm 1\} \times (\pm \{ 5,7\}) )$
can be decomposed into four $C_{16}$-factors by Theorem~\ref{d}.
$ \mathrm{Cay}(Z_m \times Z_{16},  \{ \pm 1\}  \times  \{ \pm 1, 8\} ) $
can be decomposed into two $C_{16}$-factors and a $C_{m}$-factor by Lemma~\ref{d=41}.
The last three $C_m$-factors come from a $C_m$-factorization of $ \mathrm{Cay}(Z_m \times Z_{16},  \{ \pm 1\}  \times  \{ 0, \pm 3\} )$ by Theorem~\ref{d=0d}.

\vspace{5pt}

{\bf Case 12}: $r=21$.

The 1-factor is $mI_{15}$.
$ \mathrm{Cay}(Z_m \times Z_{16},  \{ \pm 1\}  \times (\pm \{ 2,4,6\} ) ) \cup  m(\cup_{i=1}^{12} I_i)$
can be decomposed into 12 $C_{16}$-factors by Lemmas~\ref{d=2-13} and~\ref{d=6-27}.
The graph $\ m(I_{13} \cup I_{14})$ can be decomposed into a $C_{16}$-factor since $I_{13} \cup I_{14}$ can form a 16-cycle.
$ \mathrm{Cay}(Z_m \times Z_{16},  \{ \pm 1\} \times (\pm \{ 3,5,7\}) )$
can be decomposed into six $C_{16}$-factors by Theorem~\ref{d}.
$ \mathrm{Cay}(Z_m \times Z_{16},  \{ \pm 1\}  \times  \{ \pm 1, 8\} ) $
can be decomposed into two $C_{16}$-factors and a $C_{m}$-factor by Lemma~\ref{d=41}.
The last $C_m$-factor comes from a $C_m$-factorization of $ \mathrm{Cay}(Z_m \times Z_{16},  \{ \pm 1\}  \times  \{ 0\} )$ by Theorem~\ref{d=0,4}.
\qed


\begin{thebibliography}{Z}
\baselineskip 11pt

\bibitem{ABBE}
P. Adams, E. J. Billington, D. E. Bryant, and S.I. El-Zanati, On the Hamilton-Waterloo problem,
 {\it Graphs Combin.} {\bf 18} $(2002)$, 31-51.


\bibitem{ASSW}
B. Alspach,  P. J. Schellenberg,  D. R. Stinson, and D. Wagner, The
Oberwolfach problem and factors of uniform odd length cycles, {\it
J. Combin.  Theory Ser.  A} {\bf 52} $(1989)$,  20-43.


\bibitem{AKK}
J. Asplund, D. Kamin, M. Keranen, A. Pastine, and S. $\ddot{\rm O}$zkan, On the  Hamilton-Waterloo problem with triangle factors and $C_{3x}$-factors,
{\it Australas. J. Combin.} {\bf 64} $(2016)$, 458-474.


%\bibitem{AH}
%A. Assaf and A. Hartman, Resolvable group divisible designs with block size $3$,
%{\it Discrete Math. }{\bf 77} $(1989)$,  5-20.


\bibitem{BB}
S. Bonvicini and M. Buratti, Octahedral, dicyclic and special linear solutions
of some unsolved Hamilton-Waterloo problems,
{\it Ars Math. Contemp. }{\bf 14} $(2018)$,  1-14.


%\bibitem{BuD}
%M. Buratti and P. Danziger, A cyclic solution for an infinite class of Hamilton-Waterloo problems,
% {\it Graphs and Combin.}  {\bf 32} $(2016)$, 521-531.

%\bibitem{BR}
%M. Buratti and G. Rinaldi, On sharply vertex transitive $2$-factorizations of the complete graph,
% {\it J. Combin. Theory Ser. A}  {\bf 111} $(2005)$, 245-256.


%\bibitem{BDT}
%A. Burgess, P. Danziger, and T. Traetta, On the Hamilton-Waterloo problem with odd orders,
%{\it J. Combin. Des.} {\bf 25} (2017), 258-287.


\bibitem{BDT2}
A. Burgess, P. Danziger, and T. Traetta, On the Hamilton-Waterloo problem with odd cycle lengths,
{\it J. Combin. Des.} {\bf 26} (2018), 51-83.


%\bibitem{BDT3}
%A. Burgess, P. Danziger, and T. Traetta, On the Hamilton-Waterloo problem with cycle lengths of distinct parities,
%{\it Discrete Math.} {\bf 341} $(2018)$, 1636-1644.


\bibitem{BDT4}
A. Burgess, P. Danziger, and T. Traetta, The Hamilton-Waterloo problem with even cycle lengths,
{\it Discrete Math.} {\bf 342} $(2019)$, 2213-2222.


\bibitem{CNT}
H. Cao,  M. Niu, and C. Tang, On the existence of cycle frames and
almost resolvable cycle systems, {\it Discrete Math.} {\bf
311} (2011), 2220-2232.


\bibitem{DQS}
P. Danziger, G. Quattrocchi, and B. Stevens, The Hamilton-Waterloo problem for cycle sizes $3$ and $4$,
 {\it J. Combin. Des.} {\bf 17} $(2009)$, 342-352.


%\bibitem{DL2}
%J. H. Dinitz and A. C. H. Ling, The Hamilton-Waterloo problem: The case of triangle-factors and one Hamilton cycle,
% {\it J. Combin. Des.} {\bf 17} $(2009)$, 160-176.


\bibitem{FH}
H. Fu and K. Huang, The Hamilton-Waterloo problem for two even cycles factors,
 {\it Taiwanese J. Math.} {\bf 12} $(2008)$, 933-940.


\bibitem{HS}
D. G. Hoffman and P. J. Schellenberg, The existence of $C_{k}$-factorizations of $K_{2n}-F$,
 {\it Discrete Math.} {\bf 97} $(1991)$, 243-250.


% \bibitem{HNR}
% P. Horak, R. Nedela, and A. Rosa, The Hamilton-Waterloo problem: The case of Hamilton cycles and triangle-factors,
% {\it Discrete Math.} {\bf 284} $(2004)$, 181-188.


%\bibitem{K}
%D. C. Kamin, Hamilton-Waterloo problem with triangle and $C_9$-factors,
% Master's Thesis Michigan Technological University, 2011.
% https://digitalcommons/etds/207.


\bibitem{KO}
M. Keranen and S. $\ddot{{\rm O}}$zkan,  The Hamilton-Waterloo problem with 4-cycles and a singer factor of $n$-cycles,
{\it Graphs Combin.} {\bf 29} (2013), 1827-1837.

%\bibitem{KP}
%M. Keranen and A. Pastine, A generalization of the Hamilton-Waterloo problem on complete equipartite graphs,
% {\it J. Combin. Des.} {\bf 25} $(2017)$, 431-468.


\bibitem{LF}
H. Lei and H. Fu, The Hamilton-Waterloo problem for triangle-factors and heptagon-factors,
 {\it Graphs Combin.}  {\bf 32} $(2016)$, 271-278.

%  \bibitem{LFS}
% H. Lei, H. Fu, and H. Shen, The Hamilton-Waterloo problem for Hamilton cycles and $C_{4k}$-factors,
% {\it Ars. Combin.}  {\bf 100} $(2011)$, 341-348.


% \bibitem{LS}
% H. Lei and H. Shen, The Hamilton-Waterloo problem for Hamilton cycles and triangle-factors,
% {\it J. Combin. Des.} {\bf 20} $(2012)$, 305-316.


\bibitem{LD}
A. C. H. Ling and J. H. Dinitz, The Hamilton-Waterloo problem with triangle-factors and Hamilton cycles: The case $n\equiv 3\pmod{18}$,
 {\it J. Combin. Math. Combin. Comput. } {\bf 70} $(2009)$, 143-147.


%\bibitem{L2000}
%J. Liu,  A generalization of the Oberwolfach problem and
%$C_{t}$-factorizations of complete equipartite graphs,
%{\it  J. Combin. Des.} {\bf 8} (2000), 42-49.


\bibitem{L2003}
J. Liu, The equipartite Oberwolfach problem with uniform tables,
 {\it J. Combin. Theory Ser.}  {\it A} {\bf 101} $(2003)$, 20-34.


%\bibitem{MT}
%F. Merola and T. Traetta, Infinitely many cyclic solutions to the Hamilton-Waterloo problem with odd length cycles,
%{\it Discrete Math. }{\bf 339} (2016),  2267-2283.


\bibitem{OO}
U. Odaba\c{s}{\scriptsize I} and S. $\ddot{{\rm O}}$zkan, The Hamilton-Waterloo problem with $C_{4}$ and $C_{m}$ factors,
 {\it Discrete Math.} {\bf 339} $(2016)$, 263-269.


\bibitem{PWL}
W. L. Piotrowski, The solution of the bipartite analogue of the Oberwolfach problem,
{\it Discrete Math.}  {\bf 97} (1991), 339-356.


\bibitem{R}
R. Rees,  Two new direct product-type constructions for resolvable
group-divisible designs, {\it J. Combin. Des.}  {\bf 1} (1993),
15-26.


\bibitem{WC}
L. Wang and H. Cao, A note on the Hamilton-Waterloo problem with $C_8$-factors and $C_m$-factors,
                 {\it Discrete Math.} {\bf 341} (2018), 67-73.

\bibitem{WCC}
L. Wang, F. Chen, and H. Cao, The Hamilton-Waterloo problem for $C_3$-factors and $C_n$-factors,
 {\it J. Combin. Des.}  {\bf 25} (2017), 385-418.


\bibitem{WLC}
L. Wang, S. Lu, and H. Cao, Further results on the Hamilton-Waterloo problem,
  {\it J. Combin. Des.}  {\bf 26} (2018), 27-47.


\end{thebibliography}
\end{document}